\documentclass[11pt]{amsart}
\usepackage{amscd,amssymb,amsmath,graphicx}
\newtheorem{theorem}{Theorem}[section]

\theoremstyle{definition}

\theoremstyle{remark}

\numberwithin{equation}{section}
\begin{document}
\title{A combinatorial analog of a theorem of F.J. Dyson}
\author{Pallavi Jayawant}
\address{Department of Mathematics, Bates College, Lewiston,
ME 04240, U.S.A.}
\email{pjayawan@bates.edu}
\author{Peter Wong}
\address{Department of Mathematics, Bates College, Lewiston,
ME 04240, U.S.A.}
\email{pwong@bates.edu}
\begin{abstract}
Tucker's Lemma is a combinatorial analog of the Borsuk-Ulam theorem and the case $n=2$ was proposed by Tucker in 1945. Numerous generalizations and applications of the Lemma have appeared since then. In 2006 Meunier proved the Lemma in its full generality in his Ph.D. thesis. There are generalizations and extensions of the Borsuk-Ulam theorem that do not yet have combinatorial analogs. In this note, we give a combinatorial analog of a result of Freeman J. Dyson and show that our result is equivalent to Dyson's theorem. As with Tucker's Lemma, we hope that this will lead to generalizations and applications and ultimately a combinatorial analog of Yang's theorem of which both Borsuk-Ulam and Dyson are special cases.
\end{abstract}
\date{\today}
\maketitle

\newcommand{\af}{\alpha}
\newcommand{\et}{\eta}
\newcommand{\ga}{\gamma}
\newcommand{\ta}{\tau}
\newcommand{\ph}{\varphi}

\newcommand{\bt}{\beta}

\newcommand{\lb}{\lambda}

\newcommand{\wh}{\widehat}

\newcommand{\sg}{\sigma}

\newcommand{\om}{\omega}

\newcommand{\cH}{\mathcal H}

\newcommand{\cF}{\mathcal F}

\newcommand{\N}{\mathcal N}

\newcommand{\R}{\mathcal R}

\newcommand{\Ga}{\Gamma}

\newcommand{\cc}{\mathcal C}

\newcommand{\bea} {\begin{eqnarray*}}

\newcommand{\beq} {\begin{equation}}

\newcommand{\bey} {\begin{eqnarray}}

\newcommand{\eea} {\end{eqnarray*}}

\newcommand{\eeq} {\end{equation}}

\newcommand{\eey} {\end{eqnarray}}

\newcommand{\ovl}{\overline}

\newcommand{\vv}{\vspace{4mm}}

\newcommand{\lra}{\longrightarrow}


\bibliographystyle{plain}

\section{Introduction} \label{S:intro}
Sperner's Lemma \cite{sperner} and Tucker's Lemma ($n=2$ in \cite{tucker} and general $n$ in \cite{lefschetz}) are well-known combinatorial analogs of two classical theorems in topology, namely, the Brouwer Fixed Point Theorem and the Borsuk-Ulam Theorem, respectively. Fan gives a generalization of Tucker's Lemma in \cite{fan}. These lemmas have useful applications, some of which can be found in \cite{freund}, \cite{matousek}, and \cite{simmons}. In 2006, Meunier \cite{meunier} proved Tucker's Lemma in its full generality in his Ph.D. thesis. Recent work of Ziegler \cite{ziegler}, of Matou\v{s}ek \cite{matousek2} and of de Longueville and \v{Z}ivaljevi\'c \cite{longueville-zivaljevic} give further evidence that combinatorial analogs of these topological theorems are desirable as they lead to elementary and constructive proofs of these theorems and may produce algorithms that have useful applications. There are numerous generalizations and extensions of the Borsuk-Ulam theorem that do not yet have combinatorial analogs.
In this paper, we prove a combinatorial analog of the following theorem of Dyson \cite{dyson} that is in the same vein as Tucker's Lemma.

\begin{theorem} \label{yangthm} For any continuous function $f$ from $\mathbb{S}^{2}$ to $\mathbb{R}$, there exist two mutually orthogonal diameters whose four endpoints are mapped to the same value under $f$.
\end{theorem}

Our proof is constructive in the sense of Freund and Todd\cite{freund} and Prescott and Su\cite{prescott}. Furthermore, we show that our result is equivalent to Dyson's theorem.

For the combinatorial analog, we use the following terminology. A triangulation of $\mathbb{S}^2$ is $\it{symmetric}$ if for each simplex $\sigma$ in the triangulation, $-\sigma$ is also a simplex in the triangulation. In particular for the antipodal map $A$ on $\mathbb{S}^2$, $A(v)=-v$ is a vertex of the triangulation iff $v$ is a vertex. In this paper, we deal only with triangulations that are finite and symmetric.
A {\it generalized Tucker labelling of a symmetric triangulation $T$ of $\mathbb{S}^2$ by $\{\pm 1\}$} is a labelling $\ell$ of the vertices of a symmetric triangulation of $\mathbb{S}^2$ such that $\ell(-v)=-\ell(v)$ and $\ell(v)\in \{\pm 1\}$ for all vertices $v$. If the labels of two vertices sum to zero, we say that the vertices have opposite labels. Thus, a vertex and its antipode have opposite labels. Note that in a Tucker labelling the number of labels used to label the vertices is twice the dimension of the sphere. Since we are using fewer labels, we have called the labelling a generalized Tucker labelling.
The labelling can always be extended linearly to give a simplicial map on the whole sphere, i.e., a simplicial map $L:\mathbb{S}^2=|T|\to \mathbb{R}$, where $|T|$ denotes the underlying space of the triangulation, $T$. Then our main theorem is the following combinatorial analog.

\begin{theorem} \label{mainthm} For any generalized Tucker labelling $\ell$ of a symmetric triangulation of $\mathbb{S}^2$ by $\{\pm 1\}$, there exists a polygonal simple closed path that is invariant under the antipodal map $A$ and is mapped to zero under the simplicial map $L$ which is the linear extension of $\ell$.
\end{theorem}

Note that the polygonal simple closed path is not a subcomplex of the triangulation and it passes through the interiors of the simplices of the triangulation, as seen in figure 1. In section 2, we give a proof of this theorem and in section 3 we establish the equivalence with Dyson's theorem. In the concluding section we discuss the connection to Yang's theorem which is a generalization of Dyson's theorem to higher dimensions.

After the first version of this paper was completed, we learned that Kulpa, et al. in \cite{kulpa} had studied certain symmetric triangulations called proper symmetric triangulations of $\mathbb{S}^2$ using combinatorial techniques similar to our techniques. However our result is valid for all symmetric triangulations of $\mathbb{S}^2$ and consequently our proof is different from that of \cite{kulpa} because it does not use the special property of symmetric triangulations that \cite{kulpa} uses. While \cite{kulpa} establishes the existence of a maximal chain of triangles invariant under the antipodal map we construct a closed invariant path which then helps us to prove the equivalence of our theorem to Dyson's theorem. Although Meunier's idea (\cite{meunier}) could be applied to the result in \cite{kulpa} to extend it to all symmetric triangulations, more work would be needed to establish the path we construct that leads to the equivalence with Dyson's theorem. Our result is a direct approach that yields equivalence to Dyson's theorem. Note that in our proof we also obtain a chain of triangles invariant under the antipodal map.

\section{Combinatorial Proof of the Main Theorem} \label{S:proof}

We start with a generalized Tucker labelling $\ell$ of $\mathbb{S}^2$ by $\{\pm 1\}$. For any vertex $v_0$, there exists a polygonal path joining $v_0$ and $-v_0$ that goes along the edges of the triangulation. Clearly along this path there is an edge with opposite labels and consequently a triangle such that all its vertices do not have the same label.

\subsection{Existence of polygonal simple closed path} \label{S:existence}

Choose a triangle $T_1$ such that all its vertices do not have the same label. Then we have a unique line segment that is mapped to zero under the simplicial map $L$ and whose two end points $p_1$ and $p_2$ are the midpoints of two edges of $T_1$. Without loss of generality, we pick $p_1$ as the starting point of the polygonal simple closed path. We then continue the path at the other end point $p_2$ by adjoining the line segment $p_2p_3$ which is the unique line segment in the adjacent triangle $T_2$ that is mapped to zero under the simplicial map $L$ (see figure 1).

\begin{figure}[h]
\begin{center}
\includegraphics{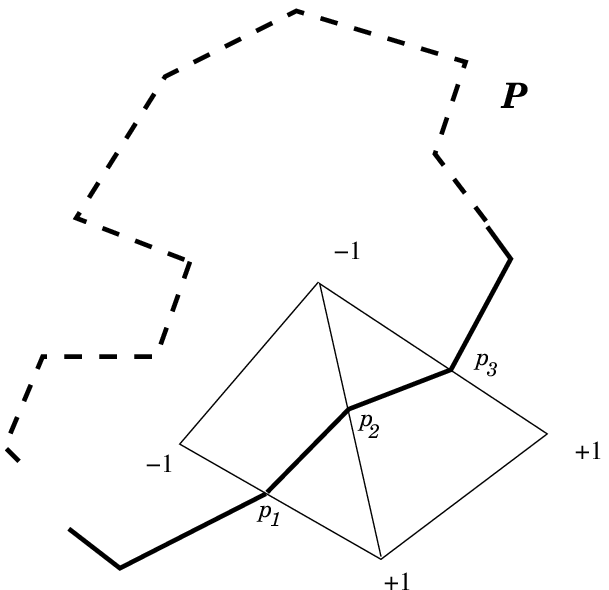}
\caption{}
\end{center}
\end{figure}

We continue to build the path by exiting one triangle and entering another triangle via the midpoint of an edge with oppositely labeled end vertices. The uniqueness of the zeros ensures that we cannot return to any triangle along the path except for the first triangle, $T_1$. Thus the construction ends to produce a simple, closed path $P$ that passes through the midpoints of the edges of a chain of triangles. By construction all the vertices in this chain of triangles on one side of $P$ have the same label and all the vertices on the other side have the opposite label.

{\bf Remark.} We note here that the existence of the triangle $T_1$ also follows from Dyson's theorem. Since $L(-x)=-L(x)$ for every $x$ in $\mathbb{S}^2$ and $L$ is a continuous real valued function on $\mathbb{S}^2$, by Dyson's theorem there exists a pair of mutually orthogonal diameters with endpoints $x, -x, y,$ and $-y$ such that $L(x)=L(-x)=L(y)=L(-y)$. So we get $L(x)=L(-x)=0$. Since $L$ maps each vertex of the triangulation to $+1$ or $-1$, the point $x$ is not a vertex of the triangulation but it belongs to a triangle $\sigma$ with all vertices not labelled the same.


\subsection{Invariance under the antipodal map} \label{S:invariance}

As we build the path $P$, if we reach a point that is the antipode of a point already on $P$, then the symmetric nature of the triangulation guarantees that the rest of the path is the antipodal image of the existing path. So in this case, the path is invariant under the antipodal map. Otherwise the path $P$ and its antipodal image $-P$ are disjoint and are contained in $L^{-1}(0)$. Since the triangulation is finite and the zeros of $L$ can only occur on a line segment joining the midpoints of two edges of a triangle, the set $L^{-1}(0)$ consists of finitely many disjoint paths $P_1$, $P_2$, \ldots, $P_m$, each path being a simple closed path which could be constructed in the same way as we constructed $P$. Our claim is that $m$ is odd so that exactly one of the paths $P_i$ is invariant under the antipodal map. Suppose $m=2k$ for some positive integer $k$. Then we can rename the paths in $L^{-1}(0)$ such that $P_{k+i}=-P_i$ for $i=1,2, \ldots, k$. By the Jordan curve theorem, the path $P_1$ separates $\mathbb{S}^2$ into two connected components, say $C_1$ and $C_2$. Note that since the path $P_1$ does not pass through the vertices of the triangulation, the triangles through which it passes get split such that one vertex of each triangle belongs to one connected component and the other two vertices belong to the other connected component (see figure 1). Now since $P_1$ and $-P_1$ are disjoint, the path $-P_1$ lies entirely in one of these components. Without loss of generality assume that $-P_1$ is contained in $C_1$. Then by the Jordan curve theorem $-P_1$ separates $C_1 \cup C_2$ (which is $\mathbb{S}^2$) into two connected components. This means $C_1$ gets split into two connected components, one of which is the antipodal image of $C_2$ and the other is invariant under the antipodal map. As with $P_1$, each triangle through which $-P_1$ passes gets splits between the two components and hence each component contains at least one vertex of the triangulation. Thus the paths $P_1$ and $-P_1$ separate 
$\mathbb{S}^2$ into three connected components, of which exactly one component is invariant under the antipodal map and the other two components are antipodal images of each other. Also each connected component contains at least one vertex of the triangulation. Now the paths $P_2$ and $-P_2$ separate the sphere further into connected components. We get five connected components of which exactly one is invariant under the antipodal map and the other four are paired up with the two components in a pair being the antipodal images of each other. Continuing with the paths $P_3$,\ldots, $P_k$ and their antipodal images, in the end we have $2k+1$ connected components of which exactly one is invariant under the antipodal map. This connected component contains a vertex and its antipode and so the simplicial map $L$ assumes both positive and negative values in this component and hence has a zero in this component. This contradicts the fact that none of the
 paths $P_i$ are in this component.


\section{Equivalence to Dyson's theorem} \label{S:equivalence}

We first prove that Theorem \ref{mainthm} is equivalent to the following theorem.

\begin{theorem} \label{invariantsubset} For any continuous function $g$ from $\mathbb{S}^{2}$ to $\mathbb{R}$ such that $g(-x)=-g(x)$, there exists a nonempty compact connected subset $X$ in $g^{-1}(0)$ that is invariant under the antipodal map $A$.
\end{theorem}

Then in \ref{S:dysonfromX} we show that Theorem \ref{invariantsubset} implies Dyson's theorem. Note that by the remark at the end of section \ref{S:existence}, Theorem \ref{mainthm} follows from Dyson's theorem because the starting point for the path in Theorem \ref{mainthm} can also be obtained by the use of Dyson's theorem. Thus we have the equivalence of Theorem \ref{mainthm} with Dyson's theorem.

\subsection{Theorem \ref{invariantsubset} $\Rightarrow$ Theorem \ref{mainthm}} \label{S:one-way}

Let $L$ be the simplicial extension of a generalized Tucker labelling $\ell$ of a symmetric triangulation of $\mathbb{S}^2$ by $\{\pm 1\}$. Then $L(-x)=-L(x)$ for every $x$ in $\mathbb{S}^2$. Also $L$ is a continuous real valued function on $\mathbb{S}^2$ and hence by Theorem \ref{invariantsubset} there exists a nonempty compact connected subset $X$ in $L^{-1}(0)$ that is invariant under the antipodal map. The definitions of $L$ and the labelling $\ell$ imply that $L^{-1}(0)$ consists of finitely many simple closed polygonal paths and since $X$ is compact, connected and invariant under the antipodal map, it has to be exactly one of these paths.

\subsection{Theorem \ref{mainthm} $\Rightarrow$ Theorem \ref{invariantsubset}} \label{S:converse}

We show the existence of a nonempty compact connected invariant subset $X$ in $g^{-1}(0)$ first in the case of a smooth function $g$ with $0$ as a regular value and then in the case of a continuous function $g$ (that is not necessarily smooth).

\subsubsection{$g$ is a smooth map with $0$ as a regular value}\label{S:P}

In this case, $g^{-1}(0)$ is a smooth one dimensional submanifold of $\mathbb{S}^2$ with no boundary. This means $g^{-1}(0)$ is the union of disjoint simple closed paths. The compactness of $g^{-1}(0)$ implies that we have finitely many components in the set and since $g(-x)=-g(x)$ the set is invariant under the antipodal map. We will show that one of the components is itself invariant under the antipodal map. We enclose each path in $g^{-1}(0)$ in an $\epsilon$-tubular neighborhood such that the set of tubular neighborhoods is also invariant under the antipodal map. We choose a finite symmetric triangulation of $\mathbb{S}^2$ so that the vertices of the triangulation do not lie in $g^{-1}(0)$. Furthermore since $g^{-1}(0)$ is a smooth submanifold, we can choose the mesh size of the triangulation small enough ($\le\epsilon$) so that any triangle that contains any part of a path in $g^{-1}(0)$ lies entirely in the $\epsilon$-tubular neighborhood enclosing the path. As a result the intersection of $g^{-1}(0)$ and every triangle $\sigma$ in the triangulation is either empty or a path that intersects the boundary of $\sigma$ at precisely two points on two distinct edges (see figure 2).

\begin{figure}[h]
\begin{center}
\includegraphics{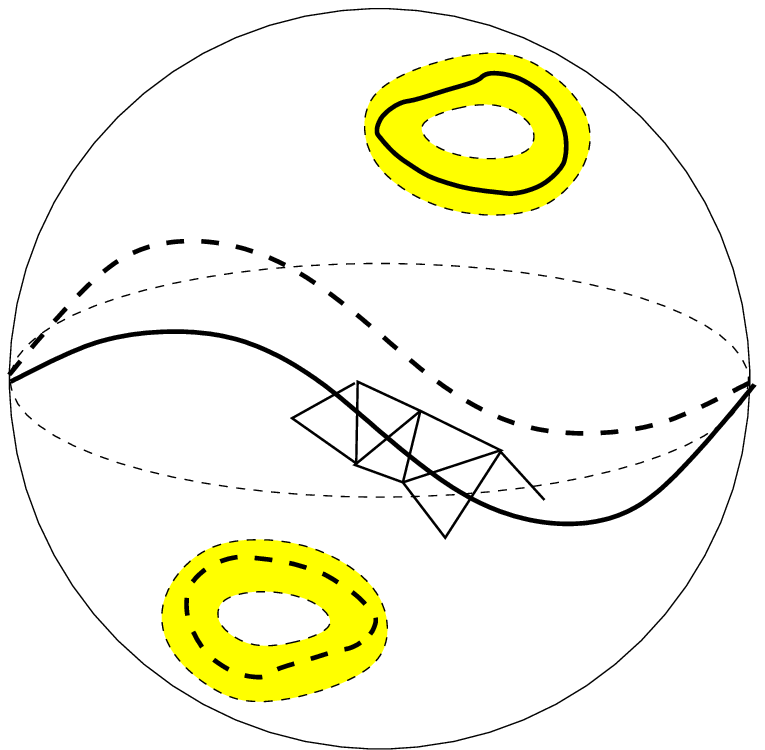}
\caption{}
\end{center}
\end{figure}

Label a vertex $v$ of the triangulation $+1$ if $g(v)$ is positive and label it $-1$ if $g(v)$ is negative. Since $g(-x)=-g(x)$, we have a generalized Tucker labelling $\ell$ of a symmetric triangulation of $\mathbb{S}^2$ by $\{\pm 1\}$. So by Theorem \ref{mainthm}, there exists a polygonal simple closed path $Q$ that is invariant under the antipodal map and is mapped to zero under the simplicial map $L$ which is the linear extension of $\ell$. From the proof of Theorem \ref{mainthm}, we know that $Q$ passes through the midpoints of some of the edges of the triangulation and the end vertices of each of these edges have opposite labels. This implies that the values of $g$ are opposite in sign at these end vertices and hence $g$ has a zero along each of these edges.
If $\sigma_1, \sigma_2, \ldots, \sigma_n, \sigma_{n+1}=\sigma_1$ is the chain of triangles through which $Q$ passes, then $g$ has at least two zeros in each $\sigma_i$. Note that from the proof of Theorem \ref{mainthm}, we know that the chain is invariant under the antipodal map.
Since the closed paths in $g^{-1}(0)$ are disjoint and the mesh size of the triangulation is sufficiently small, we have a unique simple closed path $P$ in $g^{-1}(0)$ that intersects each $\sigma_i$. Since $Q$ is invariant under the antipodal map, the path $P$ is invariant under the antipodal map. Also $P$ is compact and connected and thus we have found the required subset of $L^{-1}(0)$.

\subsubsection{$g$ is a continuous function not necessarily smooth}\label{S:X}

For every positive integer $n$, there exists a smooth real valued function $g_n$ on $\mathbb{S}^2$ (with $0$ as a regular value) such that
$|g_n(x)-g(x)|<1/n \text{ for every } x \text{ in } \mathbb {S}^2$ (\cite{hirsch}, chapter 2).
As explained in the previous paragraph, for each $g_n$ we use Theorem \ref{mainthm} and find a simple closed path $P_n$ in $g_n^{-1}(0)$ that is invariant under the antipodal map. For every $n$ and for every $x$ in $P_n$, $|g(x)|<1/n$. We consider the infinite sequence $\{P_n\}$ in the compact metric space $2^{\mathbb{S}^2}$ which is the collection of closed subsets of $\mathbb{S}^2$ with the metric given by the usual distance between two sets, i.e., dist$(A, B)=$ sup$(\{d(x,B)|x\in A\} \cup \{d(A, y)|y\in B\})$ where $A$ and $B$ are closed subsets of $\mathbb{S}^2$ and $d(x,B)=\text{min}\{\rho(x,y)|y\in B\}$, and $\rho$ is the usual metric on $\mathbb R^3$ restricted to $\mathbb{S}^2$. Let $X$ be the limit of a convergent subsequence $\{P_{n_k}\}$ of $\{P_n\}$. Since each $P_n$ is a nonempty compact connected subset of $\mathbb{S}^2$, the set $X$ is a nonempty compact connected subset of $\mathbb{S}^2$ (\cite{kuratowski}, Theorem 14). Furthermore, $X$ is invariant under the antipodal map because each $P_{n_k}$ is invariant under the antipodal map. We claim that $g(X)=0$. For each $x$ in $X$ and for every path $P_{n_k}$, there exists a point $x_{n_k}$ in $P_{n_k}$ such that dist$(P_{n_k}, x)=\rho(x_{n_k}, x)$ and hence $\rho(x_{n_k}, x)\le \text{dist}(P_{n_k},X)$. Then we have $\displaystyle\lim_{k\to\infty} |g(x_{n_k})-g(x)|=0$ because $\displaystyle\lim_{k\to\infty}\text{dist}(P_{n_k},X)=0$ and $g$ is a continuous function. But $|g(x_{n_k})|<1/{n_k}$ and hence $g(x)=0$. Thus $g(x)=0$ for each $x$ in $X$.

\subsection{Theorem \ref{invariantsubset} $\Rightarrow$ Dyson's theorem}\label{S:dysonfromX}

Let $f$ be a continuous real valued function on $\mathbb{S}^2$. Define $g:\mathbb{S}^2\to \mathbb{R}$ by $g(x)=f(x)-f(-x)$. The function $g$ is continuous and $g(-x)=-g(x)$ for all $x$ in $\mathbb{S}^2$. Then by Theorem \ref{invariantsubset} we have a nonempty compact connected subset $X$ in $g^{-1}(0)$ that is invariant under the antipodal map. Thus $f(x)=f(-x)$ for every $x$ in $X$. The rest of the proof is based on the proof of Lemma 5.5 in \cite{yang}. Since $X$ is compact, there exist points $a$ and $b$ such that $f(a)= \text{ sup } \{f(x)|x \in X\}$ and $f(b)= \text{ inf } \{f(x)|x \in X\}$. For every positive integer $n$, there exists an open covering $\mathcal{U}_n$ of $X$ such that the covering is invariant under the antipodal map and whenever $x$, $x'$ are in one set of $\mathcal{U}_n$ and $y$, $y'$ are in another set of $\mathcal{U}_n$, then $|f(x)-f(x')|<1/(2n)$ and $|\phi(x,y)-\phi(x',y')|<1/n$. For $x$ and $y$ in $\mathbb{S}^2$, $\phi(x,y)$ is the angle between the ray joining the center of the sphere and $x$ and the ray joining the center of the sphere and $y$. Since $X$ is compact and connected, we can find a sequence of points $a_0=a, a_1, \ldots,a_s=b, a_{s+1},\ldots, a_t=-a$ such that every two consecutive points are in some set in $\mathcal{U}_n$. Let $a_i=-a_{i-t}$ for $i=t+1, \ldots, 2t$, i.e., the rest of the sequence is the antipodal images of the points already in the sequence. And then we repeat this sequence of points to get an infinite sequence, i.e., let $a_j=a_i$ if $j\equiv i$ mod $2t$. For $u$ in $\mathbb R$, let $i(u)$ be the greatest integer $\le$ $tu$ and let $k_u=tu-i(u)$. Define functions $r:\mathbb{R} \to \mathbb{R}$ and $h:\mathbb{R}^2 \to \mathbb{R}$ by $r(u)=(1-k_u)\,f(a_{i(u)})+k_u\,f(a_{i(u)+1})$ and $h(u,v)=(1-k_uk_v)\,\phi(a_{i(u)}, a_{i(v)})+k_uk_v\,\phi(a_{i(u)+1}, a_{i(v)+1})$. The functions $r$ and $h$ satisfy the hypotheses of Lemma 5.2 in \cite{yang} (which is stated below) and yield a pair of real numbers $u_1^n$ and $u_2^n$ such that $r(u_1^n)=r(u_2^n)$ and $h(u_1^n, u_2^n)=\pi/2$. 

(Lemma 5.2 from \cite{yang}: Let $g:\mathbb{R}\to \mathbb{R}$ and $h:\mathbb{R}^2\to \mathbb{R}$ be maps such that, whenever $u\in \mathbb{R}$, $h(u,u)=0$, $h(u,u+1)=\pi$ and $g(u+1)=g(u)$. Let $a, b$ be elements of $\mathbb{R}$ such that $a\le b$, one of $g(a)$ and $g(b)$ is sup $\{g(u):u\in \mathbb{R}\}$ and the other is inf $\{g(u):u\in \mathbb{R}\}$. Then for any $\theta$, $0<\theta < \pi$, there exists $(u_1, u_2)\in \mathbb{R}^2$ such that $a\le u_1< u_2 < u_1+1\le b+1$, $g(u_1)=g(u_2)$ and $h(u_1, u_2)=\theta$.)

Let $x_1
^n=a_{i(u_1^n)}$ and $x_2^n=a_{i(u_2^n)}$ be points in $X$. Then the definition of $r$ and $h$ and the choice of the covering $\mathcal{U}_n$ gives $|f(x_1^n)-f(x_2^n)|<1/n$ and $|\phi(x_1^n, x_2^n)-\pi/2|<1/n$. Since $X\times X$ is compact, the sequence $\{(x_1^n, x_2^n)\}$ has a limit point $(x_1, x_2)$ in $X\times X$ which satisfies $f(x_1)=f(x_2)$ and $\phi(x_1, x_2)=\pi/2$.


\section{Conclusion} \label{S:conclusion}

We have proved Theorem \ref{mainthm} and shown that it is equivalent to Dyson's theorem. In fact, the proof in \ref{S:dysonfromX} holds not only for $\pi/2$ but for any angle $\theta$ in $(0, \pi)$. Thus Theorem \ref{mainthm} is equivalent to the following result proved independently by Livesay \cite{livesay} and Zarankiewicz \cite{zarankiewicz} of which Dyson's theorem is a special case ($r=\sqrt 2$) and hence we have the equivalence of theorems \ref{mainthm}, \ref{invariantsubset} and the following result.

\begin{theorem} For any continuous function $f$ from $\mathbb{S}^{2}$ to $\mathbb{R}$ and for each real number $r$ in $(0, 2)$, there exist points $x$ and $y$ in $\mathbb{S}^{2}$ such that $\rho(x, y)=r$ and $f(x)=f(-x)=f(y)=f(-y)$.
\end{theorem}

Dyson's theorem is also a special case ($d=1, n=2$) of the following theorem of Yang \cite{yang}:

\begin{theorem} \label{yangthm-n} For any continuous function from $\mathbb{S}^{dn}$ to $\mathbb{R}^d$, there exist $n$ mutually orthogonal diameters whose $2n$ endpoints are mapped to the same point.
\end{theorem}

The case $n=1$ is the Borsuk-Ulam theorem and its combinatorial analog is Tucker's Lemma. We have given the combinatorial analog for the case $d=1, n=2$. Su proposes that we develop combinatorial analogs for other values of $d$ and $n$ and we hope to pursue other cases. As with Tucker's Lemma, we hope that our analog leads to generalizations and applications and ultimately a combinatorial analog of Yang's theorem.

{\bf Acknowledgements} We thank Francis Su for suggesting the study of combinatorial fixed point theorems and for helpful conversations concerning this paper, and G\"unter Ziegler and an anonymous referee for helpful comments on an earlier version of the paper. We would also like to thank Jan Boronski for bringing reference \cite{kulpa} to our attention.

\end{document}